\documentclass{article}
\usepackage[a4paper]{geometry}
\usepackage[utf8]{inputenc}
\usepackage{graphicx} 
\usepackage{amsthm}
\usepackage{amsfonts}
\usepackage{amssymb}
\usepackage{amsmath}
\usepackage{mathtools}
\usepackage[maxbibnames=99]{biblatex}
\usepackage{xcolor}
\usepackage{verbatim}
\usepackage[hyphenbreaks]{breakurl}

\usepackage[citecolor=blue,colorlinks]{hyperref}

\usepackage{enumerate}

\newtheorem{innercustomgeneric}{\customgenericname}
\providecommand{\customgenericname}{}
\newcommand{\newcustomtheorem}[2]{%
  \newenvironment{#1}[1]
  {%
   \renewcommand\customgenericname{#2}%
   \renewcommand\theinnercustomgeneric{##1}%
   \innercustomgeneric
  }
  {\endinnercustomgeneric}
}

\newtheorem{thm}{Theorem}[section]

\newtheorem{Lemma}[thm]{Lemma}
\newtheorem{proposition}[thm]{Proposition}

\theoremstyle{definition}
\newtheorem{Definition}[thm]{Definition}
\newtheorem{remark}[thm]{Remark}
\newtheorem{example}[thm]{Example}

\newtheorem*{thm*}{Theorem}
\newtheorem*{corollary*}{Corollary}
\newtheorem*{lemma*}{Lemma}
\newtheorem*{proposition*}{Proposition}

\newcustomtheorem{customthm}{Theorem}

\newcommand{\bb}[1]{\mathbb{#1}}

\def\R{\mathbb{R}}
\def\to{\xrightarrow[]{}}
\def\m{\mathfrak{m}}

\def\m{\mathfrak{m}}

\mathtoolsset{showonlyrefs}

\title{On the dimension of the singular set of perimeter minimizers in spaces with a two-sided Ricci curvature bound}
\author{ \footnote{Mathematical Institute, University of Oxford, Oxford, UK. E-mail address: cucinotta@maths.ox.ac.uk} Alessandro Cucinotta
\and 
\footnote{Mathematical Institute, University of Oxford, Oxford, UK. E-mail address: fiorani@maths.ox.ac.uk - Corresponding author} Francesco Fiorani}

\begin{document}

\maketitle
\begin{abstract}
\noindent
    We show that the Hausdorff dimension of the singular set of perimeter minimizers in non-collapsed Ricci limit spaces with a two-sided Ricci curvature bound is at most $N-5$, where $N$ is the dimension of the ambient space. The estimate is sharp.
\end{abstract}

\medskip

{\small


\smallskip
\noindent
\textbf{Acknowledgements.}
Both authors are grateful to Andrea Mondino for suggesting the problem and for his useful feedback, and to Daniele Semola for some useful comments on a preliminary version of the note.
The second author is supported by the EPSRC-UKRI grant “Maths DTP 2021-22”, at the University of Oxford, with reference EP/W523781/1.
}

\section{Introduction}
We consider non-collapsed two-sided Ricci limit spaces. These are pointed metric spaces $(X,d,p)$ arising as pointed Gromov-Hausdorff limits of sequences of pointed Riemannian manifolds \\
$(M^N_k,d_k,x_k)$, where $d_k$ denotes the Riemannian distance, satisfying
\begin{equation}\label{equation | definition ricci limit in the introduction}
    |\mathrm{Ric}_{M^N_k}| \leq N-1,
\quad
\text{and}
\quad 
\mathrm{vol}(B_1(x_k)) \geq v > 0 \quad
\text{for every } k.
\end{equation}
A non-exhaustive list of works where spaces satisfying the conditions above were initially studied is \cite{AndersonDoubleRicci,BKN,TianDoubleBound,And90,CC96}.
The study of metric spaces arising as Gromov-Hausdorff limits of Riemannian manifolds with a uniform lower Ricci curvature bound, generally referred to as Ricci limit spaces, was carried out in \cite{CheegerColding1, CheegerColding2, CheegerColding3}, among others. 

We are concerned with perimeter minimizing sets in non-collapsed two-sided Ricci limit spaces. The notion of perimeter in metric spaces was studied in \cite{SomeFineProp,Am02,Mi03,AmbPall,AmbrosioDiMarinoBV}, among others. In recent years, the theory of sets of finite perimeter was further studied in  \cite{1BakryEmeryAmbrosio}, \cite{SemolaGaussGreen}, \cite{SemolaCutAndPaste} (among others) in the setting of metric measure spaces with a synthetic notion of Ricci curvature lower bounds, known as \textrm{RCD} spaces. The Riemannian Curvature Dimension condition \textrm{RCD}$(K, \infty)$ was introduced in \cite{AGS2} (see
also \cite{GigliMemo2015,AGMR}) while its finite dimensional counterpart \textrm{RCD}(K,N) was formalized in \cite{GigliMemo2015}. For a thorough introduction to the topic we refer to the survey \cite{Am02} and the references therein.

Some fundamental steps towards understanding perimeter minimizing sets in the context of \textrm{RCD} spaces were carried out in \cite{WeakLaplacian}, while other properties were then investigated in \cite{FioraniMondinoSemola}, \cite{C1}, \cite{C2}. Let us point out that minimal surfaces in Riemannian manifolds are locally boundaries of locally perimeter minimizing sets. Moreover, the study of perimeter minimizers in the \textrm{RCD} setting, due to the stability properties of the minimizing condition, allows to deduce new results about classical minimal surfaces. 

One of the key advances in the study of perimeter minimizing sets in Euclidean spaces was understanding the Hausdorff dimension of their singular set. A fundamental result obtained by De Giorgi \cite{FrontiereMinime} and refined by Federer following the work of Simons shows that a perimeter minimizing set $E \subset \R^N$ is smooth outside a closed set of Hausdorff dimension at most $N-8$. The regularity of perimeter minimizing sets in \textrm{RCD} spaces was studied in \cite{WeakLaplacian}. To report here the relevant result, we introduce some notation.

A \emph{metric measure space} is a triple $(X,d,\m)$ such that $(X,d)$ is a metric space and $\m$ is a Radon measure on $X$. We denote by $\mathcal{H}^N$ the $N$-dimensional Hausdorff measure on $(X,d)$. Given an \textrm{RCD}(K,N) space $(X,d,\mathcal{H}^N)$, we denote the set of regular points of $X$ as $R(X)$ (see the comment after Definition \ref{Definition | L^1 convergence of sets}). The regular and singular sets of a locally perimeter minimizing set $E \subset X$ (see Definition \ref{Definition | local perimeter minimizer}) are denoted respectively by $R^E$ and $S^E$ (see Definition \ref{Definition| Singular set of a finite perimeter set}). We denote by $S_{k}$ and $S^E_{k}$ respectively the $k$-dimensional singular stratum of $X$ and $E$ (see Definitions \ref{Definition | singular stratum of the ambient space} and \ref{Definition | singular stratum of a perimeter minimizer}). The space $X$ is said to be without boundary if $S_{N-2} \setminus S_{N-1} = \emptyset$.

We highlight that a point in the boundary of $E$ is regular if it is a regular point for the ambient space $X$ and, in addition, the tangent  of $E$ at $x$ is a half space. 

The sharp estimate for the Hausdorff dimension of the singular set of local perimeter minimizers in non-collapsed \textrm{RCD}(K,N) spaces shown in \cite{WeakLaplacian} states the following: if 
$(X,d,\mathcal{H}^N)$ is an \textrm{RCD}(K,N) space without boundary and $E \subset X$ is a locally perimeter minimizing set, then 
\begin{equation}\label{equation | Hausdorff dimension RCD Per min}
    \mathrm{dim}_\mathcal{H}(S^E) \leq N-3.
\end{equation}

The Hausdorff dimension of the singular set of a non-collapsed \textrm{RCD}(K,N) space without boundary is at most $N-2$, as shown in \cite{MondinoNaber}. The analogous result for non-collapsed Ricci limit spaces was shown in \cite{CheegerColding1}. In the case of non-collapsed two-sided Ricci limit spaces, a stronger estimate on the Hausdorff dimension of the singular set holds. It was shown in \cite{CheegerNaberCodim4} that if $(X,d,x)$ is a non-collapsed two-sided Ricci limit space of dimension $N$, then the Hausdorff dimension of its singular set is at most $N-4$.

Comparing the estimates for $\mathrm{dim}_\mathcal{H}(S^E)$ and $\mathrm{dim}_\mathcal{H}(S(X))$ when $E \subset X$ is a locally perimeter minimizer in a non-collapsed \textrm{RCD}(K,N) space, a question that arises is whether, under the stronger assumption that $X$ is a non-collapsed two-sided Ricci limit space, it holds \\ $\mathrm{dim}_\mathcal{H}(S^E) \leq N-5$. This is the content of the main result of this note.

\begin{customthm}{1} \label{Theorem| Hausdorff dimension of the singular set P min double bound}
Let $(X,d,p)$ be a non-collapsed two-sided Ricci limit space of dimension $N$. If $E \subset X$ is a locally perimeter minimizing set, then $S^E=S^E_{N-5}$. In particular, it holds
\begin{equation}\label{equation | Hausdorff dimension of the singular set P min double bound}
\mathrm{dim}_{\mathcal{H}}(S^E) \leq N-5.
\end{equation}
\end{customthm}

Let us point out that Theorem \ref{Theorem| Hausdorff dimension of the singular set P min double bound} is sharp, as shown in Example \ref{example|Sharpness theorem 1}. \\
Before outlining the proof of Theorem \ref{Theorem| Hausdorff dimension of the singular set P min double bound}, we mention that the key step in the proof is Theorem \ref{Theorem | Bernstein}, which is a Bernstein-type theorem for cones (see the comment before Lemma \ref{Lemma | blow up of per min in cones} for a definition of a cone) over manifolds of constant curvature $1$.

\begin{customthm}{2} \label{Theorem | Bernstein}
Let $(M,g)$ be a manifold of constant sectional curvature equal to $1$ and of dimension $N \leq 6$.
Let $C(M)$ be the metric cone over $M$ and let $p$ be its tip. If $E \subset C(M)$ is a perimeter minimizing set such that $ p \in \partial E$, then $M \cong S^N$, $C(M) \cong \bb{R}^{N+1}$ and $E \subset C(M) \cong \bb{R}^{N+1}$ is a half space.
\end{customthm}

The proof of the previous theorem follows by adapting a classical result of Simons from \cite{Simons}.
Let us mention that in Theorem \ref{Theorem | Bernstein} the assumption on the dimension is sharp, as shown by Example \ref{Example|sharpness theorem 2}. 

We now outline the proof of Theorem \ref{Theorem| Hausdorff dimension of the singular set P min double bound}.
The estimate on the Hausdorff dimension \eqref{equation | Hausdorff dimension of the singular set P min double bound} follows from the stratification result \cite[Theorem 4.1]{FioraniMondinoSemola} and $S^E = S^E_{N-5}$.

We divide the proof of $S^E = S^E_{N-5}$ in two steps. We first show that $S^E = S^E_{N-4}$. We rely on the following argument: suppose by contradiction that $x \in S^E \setminus S^E_{N-4}$. Then, by \cite[Theorem 5.12]{CheegerNaberCodim4}, the tangent space to $X$ at $x$ is isometric to $\R^N$. 
Moreover, a tangent space to $E$ at $x$ is of the form $\R^{N-3} \times A$, where $A \subset \R^3$ is itself a perimeter minimizer.
Therefore, by a classical Bernstein-type theorem (see, for instance, \cite[Theorem 28.17]{Maggi}), it follows that $A$ is a half space. This provides the desired contradiction and shows $S^E = S^E_{N-4}$.

In a second step of the proof we show $S^E_{N-4} \setminus S^E_{N-5} = \emptyset$. To this end, 
we suppose by contradiction that $x \in S^E_{N-4} \setminus S^E_{N-5}$. 
Then, by \cite[Theorem 1.16]{cheeger2018rectifiability}, a tangent space of $X$ at $x$ is $\R^{N-4} \times C(S^3/\Gamma)$, where $\Gamma \subset O(4)$ is a discrete group of isometries of the sphere acting freely. 

Moreover, a tangent space to $E$ at $x$ is of the form $\R^{N-4} \times A$, where $A \subset C(S^3/\Gamma)$ is itself a perimeter minimizer.

Therefore, by Theorem \ref{Theorem | Bernstein}, we are able to conclude $\Gamma = \{id_{S^3}\}$, $C(S^3/\Gamma) \cong \R^4$, and $A \cong \R^3 \times [0,+\infty)$, contradicting the initial assumption $x \in S^E_{N-4} \setminus S^E_{N-5}$.

\section{Preliminaries}

Throughout the paper, a \emph{pointed metric space} is a triple $(X, d, x)$, where $(X, d)$ is a complete and separable metric space and $x \in X$. We write $B_r(x)$ for the open ball centered at $x\in X$ of radius $r>0$.  Under our working conditions, the closed metric balls are compact, so we assume from the beginning the metric space $(X, d)$ to be proper. We denote by $N \geq 0$ the Hausdorff dimension of $(X,d)$, and by $\mathcal{H}^N$ the corresponding Hausdorff measure.

Given a positive Borel measure $\m$ on $X$, we denote by 
$$
L^p(X) := \left \{ u: X \to \R : \int_X |u|^p \, d \m < \infty\right\}
$$
the space of $p$-integrable functions. Given a function $u:X \xrightarrow[]{}\R$, we define its local Lipschitz constant at $x \in X$ by 
$$
    \mathrm{lip}(u)(x) := \limsup_{y\xrightarrow[]{}x} \frac{|u(x)-u(y)|}{d(x,y)}\qquad \mbox{if } x \in X \mbox{ is an accumulation point},
$$
and $\mathrm{lip}(u)(x) = 0$ otherwise.
We indicate by $\mathrm{LIP}(X)$ and $\mathrm{LIP}_\mathrm{loc}(X)$ the space of Lipschitz functions, and locally Lipschitz functions respectively.
Moreover, given a set $E \subset X$, we denote by $\chi_E$ its characteristic function.

When we refer to pointed Riemannian manifolds $(M,d,x)$, we implicitly mean that $M$ is a smooth manifold equipped with a Riemannian metric which induces the Riemannian distance $d$.

\subsection{Ricci limit spaces}

The starting point of the theory of Ricci limit spaces is the notion of Gromov-Hausdorff convergence.
\begin{Definition}[Pointed Gromov-Hausdorff convergence]\label{Definition | GH convergence}
    A sequence of pointed metric spaces $(X_k, d_k, x_k)$ is said to converge in the pointed Gromov-Hausdorff topology to $(X,d,x)$ if there exists a separable metric space $(Z,d_Z)$ and isometric embeddings $i_k : X_k \to Z$ and $i : X \to Z$ such that, for any $\varepsilon >0$ and $r>0$ there exists $\Bar{k} \in \bb{N}$ such that for for every $k \geq \Bar{k}$ we have 
    $$
        i\left(B_r^X(x)\right) \subset B_\varepsilon^Z\left(i_k\left(B_{r+\varepsilon}^{X_k}(x_k)\right)\right),
    $$
    and
    $$
        i_k\left(B_r^{X_k}(x_k)\right) \subset B_\varepsilon^Z\Big(i_k\Big(B_{r+\varepsilon}^{X}(x)\Big ) \Big).
    $$
    We denote the pointed Gromov-Hausdorff convergence by $X_k \stackrel{p\mathrm{GH}}{\longrightarrow} X$ and we say that $Z$ is the space realizing the convergence.
\end{Definition}

Let us recall the definition of the spaces we work with.
\begin{Definition}[Non-collapsed Ricci limit spaces with a two-sided bound on the Ricci curvature]\label{Definition | Ricci limit Spaces}
    Let $(M_k,d_k,x_k)$ be a sequence of pointed Riemannian manifolds of fixed dimension $N \geq 2$ such that 
        \begin{equation} \label{equation | definition non-collapsed spaces | Ricci lower bound}
        \mathrm{Ric}_{M_k} \geq N - 1.
        \end{equation}
    Suppose there exists a metric space $(X,d,x)$ such that $M_k \stackrel{p\mathrm{GH}}{\longrightarrow} X$. Then we say that $X$ is a Ricci limit space.
    If condition \eqref{equation | definition non-collapsed spaces | Ricci lower bound} is strengthened to
    \begin{equation}\label{equation | two sided ricci bound}
    |\mathrm{Ric}_{M_k}| \leq N - 1,
    \end{equation}
    we say that $X$ is a two-sided Ricci limit space. Moreover, if
    \begin{equation}\label{equation | non-collapsing assumption}
    \mathrm{vol}(B_1(x_k)) \geq v > 0
    \end{equation}
    holds, we say that $X$ is a non-collapsed (two-sided, if \eqref{equation | two sided ricci bound} holds) Ricci limit space.
\end{Definition}

Ricci limit spaces were studied extensively in \cite{CheegerColding1}, \cite{CheegerColding2} and \cite{CheegerColding3}. In particular, it was shown that the non-collapsing assumption \eqref{equation | non-collapsing assumption} forces the Hausdorff dimension of the limit space and that of the approximating sequence to be the same.

Let us recall the definition of tangent space in the setting of metric spaces.

\begin{Definition}[Tangent space to a metric space at a point] \label{Definition | tangent to Ricci limit}
    Let $(X,d)$ be a metric space and let $x \in X$. We define the \emph{space of tangent spaces} at $x$, denoted by $\mathrm{Tan}(X,x)$, to be the set of all $(Y,d_Y,y)$ such that there exists a sequence $1 > r_k > 0$, $r_k \searrow 0$ that satisfies 
    $$
    (X, d / r_k, x) \stackrel{p\mathrm{GH}}{\longrightarrow} (Y, d_Y, y).
    $$
\end{Definition}

In the case of Ricci limit spaces, Gromov compactness Theorem implies that the set of tangent spaces is always non-empty. Moreover, in \cite{CheegerColding1} it was shown that for non-collapsed Ricci limit spaces all tangent spaces are metric cones (see \cite{Burago} for the definition of metric cones).

 We now report the notion of singular stratum of a Ricci limit space.

\begin{Definition}[Singular Stratum of a Ricci limit space]\label{Definition | singular stratum of the ambient space}
    Let $X$ be a Ricci limit space. For every $k \in \bb{N}$ the $k$-singular stratum is defined to be 
    \begin{equation*}
        S_{k} := \left \{x \in X : \ \mbox{no tangent space is isometric to }\bb{R}^{k+1} \times Y \mbox{ for some metric space }Y\right\}.
    \end{equation*}
\end{Definition}

It was proved in \cite{CheegerColding1} that the Hausdorff dimension of $S^{k}$ is always less or equal to $k$.

The structure of singular strata of non-collapsed two-sided Ricci limit spaces was studied in \cite{CheegerNaberCodim4}. A key result used in the proof of Theorem \ref{Theorem| Hausdorff dimension of the singular set P min double bound} is Theorem \ref{Theorem | strati per two sided} below, which is an immediate consequence of \cite[Theorem 5.12]{CheegerNaberCodim4}.
In the next statement, if $k \in \bb{Z} \setminus \bb{N}$, we set $S_k:=\emptyset$.

\begin{thm} \label{Theorem | strati per two sided}
    Let $(X, d, x)$ be a non-collapsed two-sided Ricci limit space of dimension $N \geq 2$. Then $S_{N-1} \setminus S_{N-4}= \emptyset$.
\end{thm}

Another key result used in the proof of Theorem \ref{Theorem| Hausdorff dimension of the singular set P min double bound} is \cite[Theorem 1.16]{cheeger2018rectifiability}. More specifically, by following verbatim the proof of such theorem, one is able to conclude the following.
\begin{thm}\label{Theorem | characterization of 3 dimensional cross section}
    Let $(X, d, x)$ be a non-collapsed two-sided Ricci limit space. Then, for any $x \in S_{N-4} \setminus S_{N-5}$, there exists a tangent space at $x$ which is isometric to $\bb{R}^{N-4} \times C(S^3/\Gamma)$, where $\Gamma \subset O(4)$ is a discrete group acting freely.
\end{thm}

\subsection{Finite perimeter sets and perimeter minimizers}
A \emph{metric measure space} is a triple $(X,d,\m)$, where $(X,d)$ is a complete and separable metric space and $\m$ is a non negative Borel measure on $X$, which is finite on metric balls. Given $x \in X$, quadruples of the form $(X,d,\m,x)$ are called \emph{pointed metric measure spaces}. We say that two pointed metric measure spaces $(X,d,\m,x)$ and $(Y,\rho, \mu,y)$ are isometric in the sense of pointed metric measure spaces if there exists an isometry in the sense of metric spaces $i: X \to Y$ such that $i_{\#}\m = \mu$ and $i(x) = y$. In the case of non-collapsed Ricci limit spaces, the reference measure $\m$ corresponds to the Hausdorff measure $\mathcal{H}^N$.

For the purposes of this work, we need to introduce some background regarding finite perimeter sets and perimeter minimizers in Ricci limit spaces. In recent years, this theory has been developed in the more general framework of metric measure spaces with a synthetic notion of lower Ricci curvature bounds (RCD), which includes Ricci limit spaces. 
For an account of the theory of \textrm{RCD} spaces we refer to the survey \cite{Am02}.
 
In this section, we restrict our attention to non-collapsed \textrm{RCD}(K,N) spaces, that is \textrm{RCD}(K,N) spaces where the reference measure is the Hausdorff measure $\mathcal{H}^N$ (see \cite{GigliNonCollapsed}). The reason for this choice is that some of the results we report here were obtained in such setting. Let us remark that the family of non-collapsed Ricci limit spaces is a subset of the class of non-collapsed \textrm{RCD}(K,N) spaces.

The theory of finite perimeter sets in \textrm{RCD} spaces was developed in \cite{1BakryEmeryAmbrosio, SemolaCutAndPaste, SemolaGaussGreen, BrenaGigli, WeakLaplacian, rankOneTheorem}, among others.

We begin by recalling the definition of finite perimeter sets in the setting of metric measure spaces (see \cite{Am02,Mi03,AmbrosioDiMarinoBV}).

\begin{Definition}[Sets of locally finite perimeter]\label{definition | set of locally finite perimeter}
Let $(X,d,\m)$ be a metric measure space and let $E \subset X$ be a Borel set. Given an open set $A \subset X$, the perimeter of $E$ in $A$ is defined to be
$$
P(E, A) := \inf \left\{ \liminf_{k \to \infty} \int_A \mathrm{lip} f_k\, d\m \ : \ f_k \in \mathrm{LIP}_\mathrm{loc}(A)
, \ f_k \to \chi_E\ \mbox{ in } L^1_\mathrm{loc}(A)\right\}.
$$
The set $E \subset X$ is said to have locally finite perimeter if $P(E, B_r(x)) < + \infty$ for all $x \in X$ and $r > 0$.

Let us point out that, if $E$ has locally finite perimeter, there exists a unique Radon measure $\mu$ such that $\mu(A)=P(E,A)$ if $A \subset X$ is open. This measure is denoted $P(E,\cdot)$.
\end{Definition}

The following lemma, which follows immediately from Definition \ref{definition | set of locally finite perimeter}, is used in the proof of Theorem \ref{Theorem | Bernstein}.

\begin{Lemma} \label{Lemma | perimeters in different spaces}
    Consider two metric spaces $(X,d_x)$ and $(Y,d_y)$ with a bijective isometry $f:X \to Y$. If both spaces are equipped with Hausdorff measures of the same dimension, then for every $A,B \subset X$ we have
    \[
    P(A,B)=P(f(A),f(B)).
    \]
\end{Lemma}

In order to recall the notion of $L^1$ convergence of sets along a converging sequence of spaces introduced in \cite{AmbrosioHonda1} and \cite{1BakryEmeryAmbrosio}, we first need to report the definition pointed measured Gromov-Hausdorff convergence, which extends Definition \ref{Definition | GH convergence}.

\begin{Definition}[Pointed measured Gromov-Hausdorff convergence]\label{Definition | PMGH convergence}
    A sequence of pointed metric measure spaces $(X_k, d_k, \m_k, x_k)$ is said to converge in the pointed measured Gromov-Hausdorff topology to $(X,d,\m,x)$ if $X_k \stackrel{p\mathrm{GH}}{\longrightarrow} X$ and, using the same notation of Definition \ref{Definition | GH convergence}, we have that $({i_k})_\# \m_k \rightharpoonup i_\# \m$ w.r.t. continuous functions with bounded support on $Z$.
    We denote the pointed measured Gromov-Hausdorff convergence by $X_k \stackrel{pm\mathrm{GH}}{\longrightarrow} X$ and we say that $Z$ is the space realizing the convergence.
\end{Definition}

In the case of pointed Ricci limit spaces $(X_k,d_k,x_k)$ of the same dimension $N\geq 2$ satisfying
\[
    \mathcal{H}_k^N(B_1(x_k)) > v > 0,
\]
where $\mathcal{H}^N_k$ is the Hausdorff measure relative to $d_k$, the pointed Gromov-Hausdorff convergence (Definition \ref{Definition | GH convergence}) to a metric space $(X,d,x)$ is equivalent to the pointed measured Gromov-Hausdorff convergence (Definition \ref{Definition | PMGH convergence}) of $(X_k,d_k, \mathcal{H}_k^N, x_k)$ to $(X, d, \mathcal{H}^N, x)$. For a proof of this fact, see \cite[Theorem 1.2]{GigliNonCollapsed}.

Let us report the definition of $L^1$ convergence of sets.  

\begin{Definition}[$L^1$ convergence of sets]\label{Definition | L^1 convergence of sets}
Let $(X_k,d_k,\m_k,x_k)$ be a sequence of pointed metric measure spaces converging in pointed measured Gromov-Hausdorff sense to $(X,d,\m,x)$.
We say that a sequence of Borel sets $E_k \subset X_k$ with $\m_k(E_k) < \infty$  converges in $L^1$ sense to $E \subset X$ if, given an embedding space $(Z,d_Z)$ as in Definition \ref{Definition | PMGH convergence}, we have $\chi_{E_k}\m_k \rightharpoonup \chi_{E}\m$ in duality with continuous boundedly supported functions in $Z$ and $\m_k(E_k) \to \m(E)$. Moreover, $E_k$ is said to converge in $L^1_\mathrm{loc}$ to $E$ if the sets $E_k \cap B_r(x_k)$ converge in $L^1$ to $E \cap B_r(x)$ for all $r>0$.
\end{Definition}

Let us point out that, in the case where $(X_k,d_k, \m_k, x_k)$ is isometric to $(X, d, \m, x)$ in the sense of pointed metric measure spaces for all $k$, $L^1$ ($L^1_\mathrm{loc}$) convergence in the sense of Definition \ref{Definition | L^1 convergence of sets} is equivalent to requiring that $\m(E \Delta E_k) \to 0$ ($\m((E \Delta E_k)\cap B_r(x)) \to 0$ for all $r > 0$).

We report the definition of tangent space to a finite perimeter set found in \cite{1BakryEmeryAmbrosio}. The set of tangent spaces to a non-collapsed \textrm{RCD}(K,N) space is defined analogously to Definition \ref{Definition | tangent to Ricci limit}, where the pointed Gromov-Hausdorff convergence is replaced by convergence in the pointed measured Gromov-Hausdorff sense.

Given a non-collapsed \textrm{RCD}(K,N) space $(X,d,\mathcal{H}^N)$ and $x \in X$, we still refer to the set of tangent spaces at $x$ as \textrm{Tan}$(X,x)$. The regular set of $X$ is defined as $R(X) := \{x \in X: \mathrm{Tan}(X,x) = \{(\R^N, d_\mathrm{eucl}, 0) \}\}$. The singular set is defined as $S(X) :=X \setminus  R(X)$. By a result shown in \cite{boundary}, it holds
$$
R(X) = \{x \in X :(\R^N,d_{\mathrm{eucl}}, 0) \in \mathrm{Tan}(X,x)\}\, .
$$ 

\begin{Definition}[Tangents to a set of locally finite perimeter]\label{definition | Tangent to a set of locally finite perimeter}
Let $(X,d,\mathcal{H}^N)$ be a non-collapsed $\mathrm{RCD}$(K,N) space and let $E\subset X$ be a set of locally finite perimeter. We say that $(Y, d_Y, F, y)\in \mathrm{Tan}(X,E,x)$ if the pointed metric measure space $(Y,d_Y,\mathcal{H}^N,y)$ belongs to $\mathrm{Tan}(X,x)$ and $F\subset Y$ is a set of locally finite perimeter of positive measure such that $\chi_E$ converges in the $L^1_\mathrm{loc}$ sense of Definition \ref{Definition | L^1 convergence of sets} to $F$ along the blow-up sequence associated to the tangent $Y$.
\end{Definition}

\begin{Definition}[Perimeter minimizers]\label{Definition | local perimeter minimizer}
    Let $(X,d,\m)$ be a metric measure space and let $E \subset X$ be a set of locally finite perimeter. We say that $E$ is locally perimeter minimizing if for every $x \in X$ there exists $r>0$ with the following property: for every $F \subset X$ such that $F \Delta E \subset \subset B_r(x) $ we have $P(E,B_r(x) ) \leq P(F,B_r(x))$. \par 
    Similarly, the set $E \subset X$ is perimeter minimizing if for every $x \in X$, $r>0$ and $F \subset  X$ such that $F \Delta E \subset \subset B_r(x) $ we have that $P(E,B_r(x) ) \leq P(F,B_r(x) )$.
\end{Definition}

Let $(X,d,\mathcal{H}^N)$ be a non-collapsed \textrm{RCD}(K,N) space and suppose that $E \subset X$ is a locally perimeter minimizing set. If $(Y,d, F, y) \in \mathrm{Tan}(X, E,x)$, then $F$ is a perimeter minimizer. For a proof of this fact, we refer to \cite[Proposition 3.9]{1BakryEmeryAmbrosio}.
We report a result regarding density estimates for perimeter minimizing sets in the \textrm{RCD} setting which follows from \cite[Theorem 4.2]{DensityPerMin}.

\begin{Lemma} \label{Lemma | density estimates for perimeter minimizing sets}
Let $(X,d,\mathcal{H}^N)$ be an $\mathrm{RCD}(0,N)$ space and let $E \subset X$ be a perimeter minimizing set. Then, up to modifying $E$ on a $\mathcal{H}^N$-negligible set, there exists $C = C(N) >0$ such that for any $x \in \partial E$ and $r > 0$ 
\begin{align*}
     \frac{\mathcal{H}^N(E \cap B_r(x))}{\mathcal{H}^N(B_r(x))} \geq C, \qquad\frac{\mathcal{H}^N(B_r(x)\setminus E)}{\mathcal{H}^N(B_r(x))} \geq C.
\end{align*}
\end{Lemma}

Moreover, it follows from \cite[Theorem 4.2]{DensityPerMin} that if a set is perimeter minimizing in an \textrm{RCD}(K,N) space, then it admits both an open and a closed representative. These representatives have the same topological boundary. Whenever we refer to the boundary of a locally perimeter minimizing set, we mean the topological boundary of its open (or closed) representative.

\begin{Definition}[Regular points of locally perimeter minimizing sets]\label{Definition| Singular set of a finite perimeter set}
 Let $(X, d,\mathcal{H}^N)$ be a non-collapsed \textrm{RCD}(K,N) space and let $E\subset X$ be a locally perimeter minimizing set in the sense of Definition \ref{Definition | local perimeter minimizer}. Given $x\in \partial E$, we say that $x$ is a regular point of $E$ if
 \begin{equation}
     \mathrm{Tan}(X,E,x)=\{(\R^N,d_{\mathrm{eucl}},\R^{N-1} \times [0,+\infty),0)\}\, .
 \end{equation}
 The set of regular points of $E$ is denoted by $R^E$. The set of singular points is defined to be $S^E := X\setminus R^E$.
\end{Definition}
By the $\varepsilon$-regularity result shown in \cite[Theorem 6.8]{WeakLaplacian}, it holds
$$
R^E = \{x \in \partial E :(\R^N,d_{\mathrm{eucl}},\R^{N-1} \times [0,+\infty),0) \in \mathrm{Tan}(X,E,x)\}\, .
$$

We here recall the notion of singular stratum of a perimeter minimizing set. 

\begin{Definition}[Singular Strata]\label{Definition | singular stratum of a perimeter minimizer}
Let $(X, d,\mathcal{H}^N)$ be an \textrm{RCD}(K,N) space, $E\subset X$ a locally perimeter minimizing set  in the sense of Definition \ref{Definition | local perimeter minimizer} and $0 \leq k \leq N-3$ an integer. The $k$-singular stratum of $E$,  $S^E_k$, is defined as
\begin{equation}
    \label{singular stratum definition}
    \begin{aligned}
    S_k^E := &\{x \in \partial E: \mbox{ no element of } \mathrm{Tan}(X,E,x) \mbox{ is of the form }(Y,d_Y,  F,y), \\
    &\; \mbox{ with }(Y,d_Y,y) \mbox{ isometric to }(Z\times\R^{k+1},d_Z \times d_\mathrm{eucl},(z,0)) \mbox{ for some pointed }(Z,d_Z,z)\\
    &\; \mbox{ and }F=G\times\R^{k+1} \mbox{ with } G\subset Z \mbox{ perimeter minimizer}\}.
    \end{aligned}
\end{equation}
\end{Definition}

 A key result that we use in the proof of Theorem \ref{Theorem| Hausdorff dimension of the singular set P min double bound} is the stratification of the singular set of locally perimeter minimizing sets, which can be found in \cite[Theorem 4.5]{FioraniMondinoSemola}. This result concerns non-collapsed \textrm{RCD}(K,N) spaces with empty boundary, i.e. such that $S_{N-1} \setminus S_{N-2}=\emptyset$ (for more on boundaries of \textrm{RCD}(K,N) spaces see, for instance, \cite{BrueNaberSemola}).
 Since Ricci limit spaces have empty boundary (as shown in \cite{CheegerColding1}), the result also applies to our setting. 

\begin{thm}[Stratification of the singular set]\label{Theorem | stratification of a perimeter minimizer}
Let $(X,d,\mathcal{H}^N)$ be a non-collapsed $\mathrm{RCD}(K,N)$ space with empty boundary, and let $E \subset X$ be a locally perimeter minimizing set. Then
\begin{equation}
    \mathrm{dim}_\mathcal{H} S_k^E \leq k \quad \mbox{for }k = 0,1,..., N-3.
\end{equation}
Moreover, $S^E\setminus S_{N-2}^E = \emptyset$.
\end{thm}

Let us recall a technical lemma. A proof of the equivalent statement in Euclidean spaces can be found in \cite[Lemma $28.13$]{Maggi}, whereas a proof of Lemma \ref{L7} can be found in \cite[End of Step $3$ in the proof of Theorem $4.4$]{FioraniMondinoSemola}.

\begin{Lemma} \label{L7}
    Let $(X, d, \mathcal{H}^N)$ be a non-collapsed $\mathrm{RCD}(K,N)$ space. Let $k \in \bb{N}$ and let $A \times \bb{R}^k \subset X \times \bb{R}^k$ be a perimeter minimizing set, then $A \subset X$ is also perimeter minimizing.
\end{Lemma}

We conclude the section by reporting a useful result regarding tangent spaces to perimeter minimizing sets in cones. Given a metric space $(X,d)$, the \emph{metric cone over $X$} is defined as the set
\[
C(X):=(X \times [0,+\infty))/{\sim}
\quad
\text{with }
(x,0) \sim (y,0) \text{ for every } x,y \in X,
\]
equipped with the cone metric (see \cite{Burago} for the definition). The point $(x, 0) \in C(X)$ is called the tip of $C(X)$.
Given a metric cone $C(X)$ and a subset $A \subset X$, the cone $C(A)$ can be canonically identified with a subset of $C(X)$, and this identification is implicitly used in the statement of the next lemma.

\begin{Lemma} \label{Lemma | blow up of per min in cones}
    Let $(C(X), d)$ be a metric cone such that $(C(X),d,\mathcal{H}^N)$ is a non-collapsed $\mathrm{RCD}(0,N)$ space. If $E \subset C(X)$ is a perimeter minimizing set whose boundary contains the tip of the cone $p$, then there 
    exists a set $A \subset X$ such that $C(A) \subset C(X)$ is a perimeter minimizing cone, whose boundary contains $p$, and such that
    $$
    (C(X), d, C(A), p) \in \mathrm{Tan}(C(X), E, p).
    $$
\end{Lemma}
The proof is a consequence of the rigidity part of the Monotonicity Formula found in \cite[Theorem 3.1]{FioraniMondinoSemola} and a density estimate for locally perimeter minimizing sets shown in \cite[Lemma 5.1]{DensityPerMin}.

\subsection{Second variation formula in Euclidean spaces}
In this section we collect some technical results on perimeter minimizing sets in Euclidean spaces that are used in the proof of Theorem \ref{Theorem | Bernstein}. 
The topic is classical and we refer to \cite{Maggi} and \cite{Giusti} for an account of the theory.
Here and throughout the paper, we adopt the notation of \cite{Giusti}.  \par 
We now recall the notion of tangential derivatives to the boundary of a smooth open set $E \subset \bb{R}^N$. Let $\nu: \partial E \to S^{N-1}$ be the outward normal vector on $\partial E$ and let $g \in C^{\infty}(\bb{R}^N)$. On points of $\partial E$ the tangential derivative of $g$ is defined as
\[
\delta g:=\nabla g - \nu (\nabla g \cdot \nu).
\]
Given any integer $1 \leq i \leq N $, the $i$-th component of the tangential derivative of $g$ is defined as
\[
\delta_i g:=\partial_i g -\nu_i(\nabla g \cdot \nu).
\]
Similarly, the tangential Laplacian of $g$ is defined as
\[
\mathcal{D}g:=\sum_{i=1}^N \delta_i \delta_i g.
\]
Both $\delta g$ and $\mathcal{D} g$ depend only on the restriction of $g$ to $\partial E$; for this reason we consider tangential derivatives and tangential Laplacians of functions that are defined only on $\partial E$, assuming implicitly that we are extending the functions smoothly to $\bb{R}^N$ before applying such operators. \par
We denote by $c^2:\partial E \to \bb{R}$ the sum of the squares of the principal curvatures of $\partial E$. One can check that $c^2=\sum_{i,j}(\delta_i \nu_j)^2$ (see \cite[Remark $10.6$]{Giusti}). 
Whenever $\partial E$ is not smooth, we assume that the aforementioned objects are defined in the largest smooth subset of $\partial E$.

Moreover, we say that a set $E \subset \bb{R}^N$ is a cone with tip $p$ if it is invariant under dilations which fix $p$. This notion is consistent with the one of metric cone previously introduced. Without loss of generality, in the remainder of this work we suppose that $p$ coincides with the origin $0 \in \R^N$.

\begin{remark} \label{R2}
    By inspecting the proof of \cite[Lemma 10.9]{Giusti} one realizes that if $E \subset \bb{R}^N$ is a cone which is both smooth and has zero mean curvature in $\bb{R}^N \setminus \{0\}$, then $c^2$ is homogeneous of degree $-2$.
\end{remark}

We now recall the second variation formula for sets with vanishing mean curvature, which can be found in \cite[Identity $(10.13)$]{Giusti}.

\begin{proposition} \label{L2}
    Let $E \subset \bb{R}^N$ be an open set such that $\partial E$ is smooth and has zero mean curvature in an open set $A \subset \subset \bb{R}^N$. Let $\nu:A \to S^{N-1}$ be  an extension of the outward unit normal of $\partial E$ to $A$, and let $\zeta \in C^{\infty}_c(A)$. Define $F_t:A \to \bb{R}^N$ by $F_t(x):=x+t\zeta(x) \nu (x)$ and set $E_t:=F_t(E)$. Then
    \[
    \Big( \frac{d^2}{dt^2}P(A,E_t)\Big)_{t=0}=
    \int_{\partial E}(|\delta \zeta|^2-c^2\zeta^2) \, d\mathcal{H}^{N-1}.
    \]
\end{proposition}

We conclude the section by recalling that tangential derivatives satisfy an integration by parts formula and that $\mathcal{D} c^2$ is well behaved on minimal sets that are invariant under dilations. The next result can be found in \cite[Lemma $10.8$]{Giusti}.

\begin{Lemma} \label{L6}
    Let $E \subset \bb{R}^N$ be such that $\partial E$ is a smooth hypersurface and let $\phi \in C^{\infty}_c(\bb{R}^N)$. Then
    \[
    \int_{\partial E} \delta_i \phi \, d\mathcal{H}^{N-1}=-\int_{\partial E} \phi \nu_i \, d\mathcal{H}^{N-1}.
    \]
\end{Lemma}

The next result can be found in \cite[Lemma $10.9$]{Giusti}. 

\begin{Lemma} \label{L5}
    Let $E \subset \bb{R}^N$ be a cone which is both smooth and has vanishing mean curvature in $\bb{R}^N \setminus \{0\}$. Then in $\partial E \cap \bb{R}^N \setminus \{0\}$ we have
    \[
    \frac{1}{2} \mathcal{D} c^2 \geq-c^4+|\delta c|^2+ \frac{2c^2}{|x|^2}.
    \]
\end{Lemma}

\section{Proof of the main results}

The goal of this section is to prove Theorem \ref{Theorem| Hausdorff dimension of the singular set P min double bound} and Theorem \ref{Theorem | Bernstein}, which we restate for the reader's convenience. The proofs will be given at the end of the section.

\begin{customthm}{1} 
Let $(X,d,p)$ be a non-collapsed two-sided Ricci limit space of dimension $N$. If $E \subset X$ is a locally perimeter minimizing set, then $S^E=S^E_{N-5}$. In particular, it holds
\begin{equation}
\mathrm{dim}_{\mathcal{H}}(S^E) \leq N-5.
\end{equation}
\end{customthm}

\begin{customthm}{2} 
Let $(M,g)$ be a manifold of constant sectional curvature equal to $1$ and of dimension $N \leq 6$.
Let $C(M)$ be the metric cone over $M$ and let $p$ be its tip. If $E \subset C(M)$ is a perimeter minimizing set such that $ p \in \partial E$, then $M \cong S^N$, $C(M) \cong \bb{R}^{N+1}$ and $E \subset C(M) \cong \bb{R}^{N+1}$ is a half space.
\end{customthm}

Let us fix some notation. In this section $N \in \bb{N}$, $N \geq 2$ and $\Gamma$ is a discrete group of isometries of $S^{N-1}$ acting freely. Moreover, $\Gamma$ induces an action on $\R^N$ given in polar coordinates by $g \cdot (\omega, r):=(g(\omega), r)$. We denote by $\pi:\R^N \to \R^N/\Gamma$ the projection on the quotient space. Since $\Gamma$ acts freely on $S^{N-1}$, it follows that it also acts freely on $\R^N\setminus \{0\}$. Consequently, $\pi_{|\R^N\setminus \{0\}}$ is a covering of $(\R^N\setminus \{0\}) /\Gamma$. Therefore, it is also a local isometry.

We say that an open set $U \subset ({\R^N\setminus \{0\}}) / \Gamma$ is a \emph{cover chart} if its preimage through $\pi$ is a finite union of disjoint open sets $\{U_i\}_{i =1}^l$ (where $l$ is the cardinality of $\Gamma$) such that $\pi_{|U_i}:U_i \to U$ is a bijective isometry for every $i$. 
Given a subset $E \subset \bb{R}^N$ and $g \in \Gamma$, we denote $g \cdot E:=\{g \cdot e : e \in E\}$. Moreover, given a subset $E \subset \R^N$ ($\R^N/\Gamma$) and $t > 0$, we define the rescaled set $E/t := \{x \in \R^N (\R^N/\Gamma) : t x \in E\}$.

\begin{Definition}[$\Gamma$-symmetric sets]
We say that a set $E \subset \bb{R}^N$ is \emph{$\Gamma$-symmetric} if for every $g \in \Gamma$ we have $g \cdot E=E$.
\end{Definition}

The next lemma shows that $\Gamma$-symmetric sets arise as preimages via $\pi$ of sets in $\R^N/ \Gamma$.

\begin{Lemma}\label{Lemma | preimages of set in the quotient are Gamma symmetric}
    If $E \subset \bb{R}^N$ is $\Gamma$-symmetric, then $\pi^{-1}(\pi(E))= E$. Conversely, if $F \subset\R^N /\Gamma$, then $\pi^{-1}(F)$ is $\Gamma$-symmetric and $\pi(\pi^{-1}(F))= F$.
    \begin{proof}
    We start by showing that if $E \subset \bb{R}^N$ is $\Gamma$-symmetric, then $\pi^{-1}(\pi(E))= E$. 
    
    Observe that $\pi^{-1}(\pi(E)) \supset E$ trivially, so that we only need the other inclusion. If $x \in \pi^{-1}(\pi(E))$, then $\pi(x)=\pi(y)$ for some $y \in E$. Therefore, there exists $g \in \Gamma$ such that $g \cdot x=y$, giving that $x \in E$ as this set is $\Gamma$-symmetric. \par 
    Let us show that, if $F \subset\R^N /\Gamma$, then $\pi^{-1}(F)$ is $\Gamma$-symmetric and $\pi(\pi^{-1}(F))= F$. 
    
    Consider $x \in \pi^{-1}(F)$. We show that for every $g \in \Gamma$ we have $g \cdot x \in \pi^{-1}(F)$. To this aim, note that $\pi(x)=\pi(g \cdot x)$ so that in particular $g \cdot x \in \pi^{-1}(x) \subset \pi^{-1}(F)$. Conversely, let $x \notin \pi^{-1}(F)$. We show that for every $g \in \Gamma$ it holds $g \cdot x \notin \pi^{-1}(F)$. Indeed, if $g \cdot x \in \pi^{-1}(F)$, then $x \in \pi^{-1}(F)$. \par 
    Finally, $\pi(\pi^{-1}(F))= F$ since $\pi$ is surjective.
    \end{proof}
\end{Lemma}

\begin{Definition}[$\Gamma$-symmetric sets minimizing the perimeter against $\Gamma$-symmetric competitors] \label{Definition | Gamma symmetric perimeter minimizers}
    We say that a $\Gamma$-symmetric set $E \subset \bb{R}^N$ minimizes the perimeter against $\Gamma$-symmetric competitors if for every $\Gamma$-symmetric set $A \subset \bb{R}^N$ and $r>0$ such that $E \Delta A \subset \subset B_r(0)$ we have
    \[
    P(E,B_r(0)) \leq P(A,B_r(0)).
    \]
\end{Definition}

The following key lemma allows us to compare the perimeter of subsets of $\R^N/\Gamma$ with the perimeter of their preimage through the projection map $\pi$.

\begin{Lemma} \label{Lemma | perimeter of Gamma symmetric}
    If $F \subset \R^N/\Gamma$ and $l \in \bb{N}$ is the cardinality of $\Gamma$, then for every measurable set $U \subset \R^N/\Gamma$, it holds
       \[
        l P(F,U)=P(\pi^{-1}(F),\pi^{-1}(U)).
       \]
       \begin{proof}
        We claim that we can find a countable collection of disjoint measurable subsets $\{B_i\}_{i \in \bb{N}}$ of $\R^N/\Gamma$ which covers $(\R^N \setminus \{0\}) / \Gamma$ and is such that each set $B_i$ is contained in a cover chart of $(\R^N \setminus \{0\}) / \Gamma$. To this aim, let $\{A_i\}_{i \in \bb{N}}$ be a covering of $(\R^N \setminus \{0\}) / \Gamma$ with cover charts, which exists as every point has a neighborhood which is a cover chart. To obtain a disjoint cover we define $B_1:=A_1$ and $B_{i+1}:=A_{i+1} \setminus \cup_{j=1}^i A_j$.
        
        For every $i \in \bb{N}$ the preimage $\pi^{-1}(A_i)$ coincides with the disjoint union $\cup_{j=1}^l A_i^j$, and for every integer $1 \leq j \leq l$
        $$
        \pi_{|A_i^j}:A_i^j \to A_i
        $$ 
        is a bijective isometry. 
        In particular, by Lemma \ref{Lemma | perimeters in different spaces} 
        \[
        P(F,B_i \cap U)=P((\pi_{|A_i^j})^{-1}(F),(\pi_{|A_i^j})^{-1}(B_i \cap U))
        \quad 
        \text{for every } j = 1,...,l.
        \]
        Since 
        \[
        A_i^j \cap (\pi_{|A_i^j})^{-1}(F)=
        A_i^j \cap \pi^{-1}(F),
        \]
        it holds
        \[
        P(F,B_i \cap U)=P(\pi^{-1}(F),(\pi_{|A_i^j})^{-1}(B_i \cap U))
        \quad 
        \text{for every } j = 1,...,l.
        \]
        Summing over $j$ we then get that for every $i \in \bb{N}$ there holds 
        $$
        l P(F,B_i \cap U)=P(\pi^{-1}(F),\pi^{-1}(U \cap B_i)).
        $$
        The collection $\{\pi^{-1}(B_i)\}_{i \in \bb{N}}$ is also a disjoint cover of $\bb{R}^N \setminus \{0\}$ so that, taking into account that $P(F,\{0\})=P(\pi^{-1}(F),\{0\})=0$, we obtain
        \[
        l P(F,U)= \sum_{i \in \bb{N}} l
        P(F,U \cap B_i)=
         \sum_{i \in \bb{N}} P(\pi^{-1}(F),\pi^{-1}(B_i \cap U))=
        P(\pi^{-1}(F),\pi^{-1}(U)).
        \]
        \end{proof}
\end{Lemma}

The next lemma shows that there exists a correspondence between perimeter minimizers in $\R^N/\Gamma$ and $\Gamma$-symmetric sets minimizing the perimeter against $\Gamma$-symmetric competitors in $\bb{R}^N$.

\begin{Lemma}\label{Lemma | equivalence of Gamma symmetric minimizers and minimizers in the quotient}
    Let $F \subset \R^N/\Gamma$ be a perimeter minimizing set, then $\pi^{-1}(F)$ is a $\Gamma$-symmetric set minimizing the perimeter against $\Gamma$-symmetric competitors in $\bb{R}^N$. 
    Conversely, if $E \subset \R^N$ is a $\Gamma$-symmetric set minimizing the perimeter against $\Gamma$-symmetric competitors, then $\pi(E) \subset \R^N / \Gamma$ is a perimeter minimizing set.
    \begin{proof}
    Let $F \subset \R^N/\Gamma$ be a perimeter minimizing set. Then $\pi^{-1}(F)$ is a $\Gamma$-symmetric set by Lemma \ref{Lemma | preimages of set in the quotient are Gamma symmetric}. We now show that $\pi^{-1}(F)$ also minimizes the perimeter with respect to $\Gamma$-symmetric competitors. Let $r > 0$ and $E' \subset \R^N$ be a $\Gamma$-symmetric set such that $\pi^{-1}(F) \Delta E' \subset \subset B_r(0)$. Since $\pi(B_r(0)) = B_r(0) \subset \R^N / \Gamma$, using Lemma \ref{Lemma | perimeter of Gamma symmetric} we obtain
    \[
    P(\pi^{-1}(F),B_r(0))=l^{-1} P(F,B_r(0)) \leq l^{-1} P(\pi(E'),B_r(0))=
    P(E',B_r(0)).
    \]
    Since $r > 0$ is arbitrary, we conclude that $\pi^{-1}(F)$ is a $\Gamma$ symmetric set minimizing the perimeter against $\Gamma$-symmetric competitors in $\R^N$.
    
    In an analogous fashion, one can show that if $E \subset \R^N$ is a $\Gamma$-symmetric set minimizing the perimeter against $\Gamma$-symmetric competitors, then $\pi(E) \subset \R^N / \Gamma$ is a perimeter minimizing set.
    \end{proof}
\end{Lemma}

The next proposition shows that $\Gamma$-symmetric sets minimizing the perimeter against $\Gamma$-symmetric competitors in $\R^N$ are locally perimeter minimizing in $\R^N \setminus \{0\}$.

\begin{proposition} \label{P1}
    Let $E \subset \R^N$ be a $\Gamma$-symmetric set minimizing the perimeter against $\Gamma$-symmetric competitors in $\R^N$, then $E$ is locally perimeter minimizing in $\bb{R}^N \setminus \{0\}$. In particular, in $\R^N \setminus \{0\}$ the set $E$ admits an open and a closed representative sharing the same topological boundary. Moreover, if $N \leq 7$, then $E$ has smooth boundary with vanishing mean curvature in $\bb{R}^N \setminus \{0\}$.
    \begin{proof}
    By Lemma \ref{Lemma | equivalence of Gamma symmetric minimizers and minimizers in the quotient}, the set $\pi(E) \subset \R^N /\Gamma$ is perimeter minimizing.
        Since the restricted projection map $\pi:\bb{R}^N \setminus \{0\} \to (\R^N/\Gamma) \setminus \{0\}$ is a local isometry, $E$ is then locally perimeter minimizing.
    \end{proof}
\end{proposition}

When we refer to a $\Gamma$-symmetric set minimizing the perimeter against $\Gamma$-symmetric competitors in $\R^N$, we implicitly mean its open representative. 
In the next lemma we deal with tangent spaces to sets of finite perimeter in $\R^N$ and $\R^N/\Gamma$. In both cases, when referring to elements of the tangent space at a point, we omit the distance.

\begin{Lemma} \label{L4}
    Let $E \subset \R^N$ be a $\Gamma$-symmetric set which minimizes the perimeter against $\Gamma$-symmetric competitors and whose boundary contains $0$.
    Then there exists a $\Gamma$-symmetric cone $E' \subset \R^N$ which minimizes the perimeter against $\Gamma$-symmetric competitors, whose boundary contains $0$, and such that
    \[
    (\R^N,E',0) \in  \mathrm{Tan}(\R^N, E, 0).
    \]
    \begin{proof}
        By Lemma \ref{Lemma | equivalence of Gamma symmetric minimizers and minimizers in the quotient}, $\pi(E) \subset \bb{R}^N/\Gamma$ is a perimeter minimizing set. In particular, by Lemma \ref{Lemma | blow up of per min in cones}, there exists a perimeter minimizing cone $\pi(E)^\infty \subset \R^N /\Gamma$, whose boundary contains $0$, and such that
        \[
        (\R^N/\Gamma,\pi(E)^\infty ,0) \in \mathrm{Tan}(\R^N/\Gamma,\pi(E) ,0).
        \]
We claim that
 \begin{equation} \label{Equation|preimage of tangent cone}
    (\R^N,\pi^{-1}(\pi(E)^\infty),0) \in  \mathrm{Tan}(\R^N, E, 0).
    \end{equation}
    We fix $x \in (\R^N \setminus \{0\})/\Gamma$ and we consider a bounded cover chart $A \subset (\R^N \setminus \{0\}) /\Gamma$ containing $x$. The preimage $\pi^{-1}(A)$ coincides with the disjoint union of open sets $\cup_{j=1}^l A_j$ such that the restricted maps $\pi_{|A_j}:A_j \to A$ are bijective isometries.
    Since 
    $(\R^N/\Gamma,\pi(E)^\infty ,0) \in \mathrm{Tan}(\R^N/\Gamma,\pi(E) ,0)$, there exists a sequence $t_k \to 0$, independent of $A$, such that
    \[
    \|1_{\pi(E)/t_k}-
    1_{\pi(E)^\infty}\|_{L^1(A)} \to 0 \quad \text{as }t_k \to 0.
    \]
    Taking into account that $\pi^{-1}(\pi(E)/t_k)=E/t_k$, for every $j=1,...,l$ it holds
    \begin{equation} \label{Equation|blow up}
    \|1_{E/t_k}-
    1_{\pi^{-1}(\pi(E)^\infty)}\|_{L^1(A_j)} \to 0 \quad \text{as }t_k \to 0.
    \end{equation}
    Hence, we can construct a locally finite open cover $\{A_j\}_{j \in \bb{N}}$ of $\R^N \setminus \{0\}$, such that for every $j \in \bb{N}$ condition \eqref{Equation|blow up} is satisfied. Setting $B_1:=A_1$ and $B_j:=A_j \setminus \cup_{i=1}^{j-1}B_i$ we obtain a refinement of the cover $\{A_j\}_{j \in \bb{N}}$ consisting of disjoint sets. Since $\{A_j\}_{j \in \bb{N}}$ is locally finite in $\R^N \setminus \{0\}$ also $\{B_j\}_{j \in \bb{N}}$ has this property.
    Hence, for every $r,\varepsilon>0$ with $r>\varepsilon$, there exists a finite subset $I_{r,\varepsilon} \subset \bb{N}$ such that $\{B_j\}_{j \in I_{r,\varepsilon}}$ covers $B_r(0) \setminus B_{\varepsilon}(0)$. Having fixed $r>\varepsilon>0$ we then obtain
    \begin{align} \label{Equation|blow up 2}
     \|1_{E/t_k}-
    1_{\pi^{-1}(\pi(E)^\infty)} &\|_{L^1(B_r(0))} \\
    &\leq \|1_{E/t_k}-
    1_{\pi^{-1}(\pi(E)^\infty)}\|_{L^1(B_\varepsilon(0))}+ \sum_{j \in I_{r,\varepsilon}}
     \|1_{E/t_k}-
    1_{\pi^{-1}(\pi(E)^\infty)}\|_{L^1(B_j)}.
    \end{align}
    Since $E/t_k$ and $\pi^{-1}(\pi(E)^\infty)$ are both perimeter minimizing, it follows from Lemma \ref{Lemma | density estimates for perimeter minimizing sets} that 
    \[
    \sup_{k \in \bb{N}} \|1_{E/t_k}-
    1_{\pi^{-1}(\pi(E)^\infty)}\|_{L^1(B_\varepsilon(0))}
     \to 0 \quad \text{as }\varepsilon \to 0.
    \]
    On the other hand, since $I_{r,\varepsilon} \subset \bb{N}$ is finite, it holds
    \[
    \sum_{j \in I_{r,\varepsilon}}
    \|1_{E/t_k}-
    1_{\pi^{-1}(\pi(E)^\infty)}\|_{L^1(B_j)}
    \to 0 \quad \text{as }t_k \to 0.
    \]
    Hence, passing to the limit in \eqref{Equation|blow up 2}, one obtains
    \[
    \|1_{E/t_k}-
    1_{\pi^{-1}(\pi(E)^\infty)}\|_{L^1(B_r(0))} \to 0 \quad \text{as }t_k \to 0,
    \]
    proving claim \eqref{Equation|preimage of tangent cone}. \par 
    Since $\pi(E)^\infty \subset \R^N/\Gamma$ is a cone, $\pi^{-1}(\pi(E)^\infty) \subset \R^N$ is also a cone. Since $0 \in \partial (\pi(E)^\infty) \subset \R^N/\Gamma $, then $0 \in \partial (\pi^{-1}(\pi(E)^\infty)) \subset \R^N$ as well. Finally, since $\pi(E)^\infty \subset \R^N/\Gamma$ is a perimeter minimizer, $\pi^{-1}(\pi(E)^\infty) \subset \R^N$ is a $\Gamma$-symmetric set minimizing the perimeter against $\Gamma$-symmetric competitors by Lemma \ref{Lemma | equivalence of Gamma symmetric minimizers and minimizers in the quotient}.
    \end{proof}
\end{Lemma}

 We recall that, given a smooth set $E \subset \R^N$, we denote by $c^2:\partial E \to \bb{R}$ the sum of the squares of the principal curvatures of $\partial E$.

\begin{Definition}[$\Gamma$-symmetric functions]
Given a $\Gamma$-symmetric set $E \subset \bb{R}^N$ we say that a function $f:E \to \bb{R}$ is $\Gamma$-symmetric if $f(x)=f(g \cdot x)$ for every $x \in E$ and $g \in \Gamma$.
\end{Definition}

\begin{Lemma} \label{L3}
    Let $E \subset \bb{R}^N$ be a $\Gamma$-symmetric cone which is smooth in $\bb{R}^N \setminus \{0\}$ and has zero mean curvature in $\bb{R}^N \setminus \{0\}$. There exists a smooth $\Gamma$-symmetric extension of $c^2:\partial E \setminus \{0\} \to \bb{R}$ to $\bb{R}^N\setminus \{0\}$.
    \begin{proof}
        Since $E$ is a cone which is smooth in $\R^N \setminus \{0\}$, the intersection $S:=\partial E \cap S^{N-1}$ is a $N-2$ dimensional closed smooth manifold in $S^{N-1}$. We now consider the function $c^2$ restricted to $S$. We show that it can be extended to a function $h:S^{N-1} \to \R$ with the property that for every $g \in \Gamma$ and every $x \in S^{N-1}$ it holds $h(g \cdot x)=h(x)$.
        
        Let $U$ be a tubular neighborhood of $S$ in $S^{N-1}$ and let $\pi_S:U \to S$ be the nearest point projection in $S^{N-1}$.  
        We define $h' \in C^\infty(U)$ by 
        $
        h'(x):=c^2(\pi_S(x)).
        $
        Let $\eta \in C^\infty_c(U)$ be a function on $U$ which depends only on the distance from $S$ and is identically equal to $1$ on $S$.
        We define $h \in C^\infty(S^{N-1})$ by
        \[
        h(x):=
        \begin{cases}
            \eta(x) h'(x) & x \in U \\
            0 & x \notin U.
        \end{cases}
        \]
        Since $S \subset S^{N-1}$ is invariant under the action of $\Gamma$ on $S^{N-1}$, so is $U$. Similarly, the restricted function $c^2 \in C^\infty(S)$ has the property that for every $g \in \Gamma$ and every $x \in S$ it holds $c^2(g \cdot x)=c^2(x)$.
        Consequently, the extension $h' \in C^\infty(U)$ is also invariant under the action of $\Gamma$ on $U$. The same is true for $h \in C^\infty(S^{N-1})$ given the choice of $\eta \in C^\infty_c(U)$.
        
        Finally, we define the extension $f \in C^\infty(\R^N\setminus \{0\})$ in polar coordinates by $f(\omega,r):=r^{-2}h(\omega)$. Since $c^2:\partial E \setminus \{0\} \to \bb{R}$ is homogeneous of degree $-2$ by Remark \ref{R2}, the function $f \in C^\infty(\R^N \setminus \{0\})$ is the desired extension.
    \end{proof}
\end{Lemma}

We here prove Theorem \ref{Theorem | Bernstein}. Using the results we have shown so far we are able to reduce ourselves to studying $\Gamma$-symmetric sets in $\R^N$ which minimize the perimeter with respect to $\Gamma$-symmetric competitors. We are then able to follow the computations found in \cite{Simons} (see, for instance, \cite[Theorem 10.10]{Giusti}) to conclude that such sets are half spaces.

\medskip

    \textit{Proof of Theorem \ref{Theorem | Bernstein}.}
        By a standard classification result regarding manifolds of constant sectional curvature (see, for instance, \cite[Theorem 4.1]{do1992riemannian}), $M$ is isometric to $S^N / \Gamma$, where $\Gamma$ is a discrete group of isometries of $S^N$ acting freely. In conformity with the rest of this section, we consider the induced action of $\Gamma$ on $\R^{N+1}$. We refer to the associated projection as $\pi: \R^{N+1} \to \R^{N+1}/\Gamma$ .
        
        Let us point out that $C(S^N/\Gamma)$ is isometric to $\R^{N+1}/\Gamma$. Therefore, to prove the statement of the theorem it is sufficient to show the following:
        if $F \subset \R^{N+1}/\Gamma$ is a perimeter minimizing set such that $0 \in \partial F$, then $\Gamma=\{id_{S^N}\}$, and $F \subset \R^{N+1}$ is a half space.
        
        By Lemma \ref{L4} there exists a $\Gamma$-symmetric cone $G \subset \R^{N+1}$ which minimizes the perimeter against $\Gamma$-symmetric competitors, whose boundary contains $0$, and such that
        \[
        (\R^{N+1},G,0) \in  \mathrm{Tan}(\R^{N+1}, \pi^{-1}(F), 0).
        \] \par 
       
        By Lemma \ref{P1} the cone $G$ is smooth with vanishing mean curvature except at $\{0\}$. We follow the computations of \cite{Simons} (see, for instance, \cite[Theorem $10.10$]{Giusti}) for perimeter minimizing cones in $\bb{R}^{N+1}$ to show that $G$ is a half space. \par 
        If $g \in C^{\infty}_c(\bb{R}^{N+1} \setminus \{0\})$ is a $\Gamma$-symmetric function, then the set $G_t$ (using the notation of Lemma \ref{L2}) is $\Gamma$-symmetric. 
        In particular, this holds if $g \in C^{\infty}_c(\bb{R}^{N+1}\setminus \{0\})$ is a radial function. 
        Let us denote by $c^2$ the $\Gamma$-symmetric extension of $c^2$ to $\bb{R}^{N+1}\setminus \{0\}$ obtained in Lemma \ref{L3}.
        Applying Lemma \ref{L2} 
        to the $\Gamma$-symmetric product $gc$, where $g \in C^{\infty}_c(\bb{R}^{N+1}\setminus \{0\})$ is a radial function, we have
        \[
        \int_{\partial G}
        |\delta (gc)|^2-c^2(gc)^2 \, d\mathcal{H}^{N} \geq 0.
        \]
        Using Lemma \ref{L6} we obtain
        \begin{align*}
        \int_{\partial G} c^4g^2 \, d\mathcal{H}^{N} &\leq 
        \int_{\partial G} \Big(c^2 |\delta g|^2+g^2 |\delta c|^2 + \frac{1}{2} \delta c^2 \cdot \delta g^2 \Big) \, d\mathcal{H}^{N} \\
        &= \int_{\partial G} \Big(c^2 |\delta g|^2+g^2 |\delta c|^2 - \frac{1}{2} g^2 \mathcal{D}c^2 \Big) \, d\mathcal{H}^{N}.
        \end{align*}
        From Lemma \ref{L5} it follows that
        \begin{equation}\label{equation | Giusti's key inequality}
        \int_{\partial G}\Big( |\delta g|^2-\frac{2g^2}{|x|^2} \Big)c^2 \, d\mathcal{H}^{N} \geq 0.         
        \end{equation}
        The same is true by approximation for every radial function $g \in C^{\infty}(\bb{R}^{N + 1})$ such that
        \begin{equation}\label{equation | Giusti bound}
            \int_{\partial G}\frac{g^2}{|x|^2}c^2 \, d\mathcal{H}^{N} < + \infty.
        \end{equation}
        Since $c^2$ is homogeneous of degree $-2$ by Remark \ref{R2}, the previous condition holds if \\ $g \in C^\infty(\R^{N+1})$ satisfies
        \begin{equation}\label{equation | Giusti bound 2}
            \int_{\partial G}\frac{g^2}{|x|^4}\, d\mathcal{H}^{N} < + \infty.
        \end{equation}
        A function $g \in C^\infty(\R^{N+1})$ of the form
        \[
        g(x):=|x|^\alpha \max\{|x|,1\}^\beta
        \]
        satisfies \eqref{equation | Giusti bound 2} if
        \begin{equation}\label{equation | system for alpha and beta}
            \begin{cases}
                \alpha > \frac{4-N}{2}, \\
                \alpha+\beta <\frac{4-N}{2}.
            \end{cases}
        \end{equation}
        Plugging such $g$ in \eqref{equation | Giusti's key inequality}
        we obtain
        \begin{equation}\label{equation | final inequality for half space}
        (\alpha^2-2)\int_{\partial G \cap B_1(0)}|x|^{2\alpha-2}c^2 \, d\mathcal{H}^N
        +
        ((\alpha+\beta)^2-2)\int_{\partial G \setminus B_1(0)} |x|^{2(\alpha+\beta)-2}c^2 \, d \mathcal{H}^N \geq 0.
        \end{equation}
        Since $N \leq 6$ we can choose $\alpha$ and $\beta$ compatible with \eqref{equation | system for alpha and beta} and such that $\alpha^2-2 \leq 0$ and $(\alpha+\beta)^2-2 \leq 0$. With such choice of $\alpha$ and $\beta$, from inequality \eqref{equation | final inequality for half space} it follows that $c^2$ is identically $0$ on $\partial G$. Hence, $G \subset \R^{N+1}$ is a half space.

        Let us show that $\Gamma = \{id_{S^{N}}\}$. Since $G$ is a half space, $\partial G \cap S^{N}$ is a great circle. Since $G$ is $\Gamma$-symmetric, $\partial G \cap S^{N}$ is sent to itself by all elements of $\Gamma$. We claim that the poles with respect to this great circle (that is, the two points on $S^{N}$ at maximal distance from $\partial G \cap S^{N}$) are swapped by every element of $\Gamma$ which is not the identity. Indeed, the poles cannot be fixed as the action of $\Gamma$ is free. Furthermore, the distance between each pole and $\partial G\cap S^{N}$ must be preserved since all elements of $\Gamma$ are isometries. In particular, all the elements of $\Gamma$ which are not the identity swap the poles. 
        \par
        Since $\Gamma$ acts on $\R^N$ isometrically, for every $g \in \Gamma$, $g \neq id_{S^N}$, there is a neighborhood $U \subset \R^N$ of one of the poles such that $U \subset G$ and $g \cdot U \subset \R^N \setminus G$. Consequently, any such element of $\Gamma$ cannot map the half space $G$ to itself. As $G$
is $\Gamma$-symmetric, we conclude that $\Gamma$ is trivial.
\par 
        Finally, the initial set $F \subset \R^{N+1}/\Gamma \cong \R^{N+1}$ is a half space since there are no nontrivial perimeter minimizing sets in Euclidean spaces of dimension less than $8$.
    \qed

\medskip

Theorem \ref{Theorem | Bernstein} fails if $N \geq 7$ as shown by the next example.

\begin{example} \label{Example|sharpness theorem 2}
Let $N=7$ and let $\Gamma:=\{id_{S^7},-id_{S^7}\}$. Let $E \subset \R^8$ be the Simons cone. Let us note that $E \subset \R^{8}$ is a $\Gamma$-symmetric set which minimizes the perimeter (and, in particular, it minimizes the perimeter against $\Gamma$-symmetric competitors). By Lemma \ref{Lemma | equivalence of Gamma symmetric minimizers and minimizers in the quotient}, $\pi(E) \subset \R^{8}/\Gamma \cong C(S^7/\Gamma)$ is a perimeter minimizing set. Moreover, the boundary of $\pi(E)$ contains $0$.

Let now $N=7+k$ with $k \in \bb{N}$
and let $\Gamma:=\{id_{S^{7+k}},-id_{S^{k+7}}\}$. Let $E \times \R^k \subset \R^8 \times \R^k$ be the product of the Simons cone with the extra Euclidean factors. $E \times \R^k \subset \R^8 \times \R^k$ is a $\Gamma$-symmetric set which minimizes the perimeter (and, in particular, it minimizes the perimeter against $\Gamma$-symmetric competitors). By Lemma \ref{Lemma | equivalence of Gamma symmetric minimizers and minimizers in the quotient}, $\pi(E) \subset \R^{8+k}/\Gamma \cong C(S^{7+k}/\Gamma)$ is a perimeter minimizing set. Moreover, the boundary of $\pi(E)$ contains $0$.
\end{example}

Building on Theorem \ref{Theorem | Bernstein}, we prove Theorem \ref{Theorem| Hausdorff dimension of the singular set P min double bound}.
\medskip

    \textit{Proof of Theorem \ref{Theorem| Hausdorff dimension of the singular set P min double bound}.}
The proof is divided in two steps: we start by showing $S^E \setminus S^E_{N-4} = \emptyset$, and then prove $S^E_{N-4} \setminus S^E_{N-5} = \emptyset$.

Let us show $S^E \setminus S^E_{N-4} = \emptyset$. Suppose by contradiction that $x \in S^E \setminus S^E_{N-4} $. Then, there exists a pointed metric space $(Y,d_y,p)$ and a perimeter minimizing set of the form $\R^{N-3} \times A \subset \R^{N-3} \times Y$ whose boundary contains $(0,p) \in \R^{N-3} \times Y$ such that 
$$
(\R^{N-3} \times Y, \R^{N-3} \times A,(0,p)) \in \mathrm{Tan}(X,E,x).
$$
By Lemma \ref{L7}, $A \subset Y$ is a perimeter minimizing set whose boundary contains $p$. Moreover, $\R^{N-3} \times Y \cong \R^N$ by Theorem \ref{Theorem | strati per two sided} since 
$$
(\R^{N-3} \times Y, (0,p)) \in \mathrm{Tan}(X,x).
$$ 
Therefore, $A \subset \R^3$ is a half space as it minimizes the perimeter and has non-empty boundary. We conclude that $x \not \in S^E \setminus S^E_{N-4}$, a contradiction.

In the rest of the proof we show
$S^E_{N-4} \setminus S^E_{N-5} = \emptyset$. 
Suppose by contradiction that there exists $x \in S^E_{N -4} \setminus S^E_{N-5}$. 

By Theorem \ref{Theorem | characterization of 3 dimensional cross section} the tangent space $\mathrm{Tan}(X,x)$ has an element isometric to $\bb{R}^{N-4} \times C(S^3/\Gamma)$, where $\Gamma \subset O(4)$ is a discrete group acting freely.
Hence, 
\[
(\bb{R}^{N-4} \times C(S^3/\Gamma), \R^{N-4} \times F,(0,p)) \in \mathrm{Tan}(X,E,x),
\]
where $p \in C(S^3/\Gamma)$ is the tip, and $\R^{N-4} \times F \subset \bb{R}^{N-4} \times C(S^3/\Gamma)$ is a perimeter minimizing set whose boundary contains $(0,p) \in \R^{N-4} \times C(S^3/\Gamma)$.  

By Lemma \ref{L7}, $F \subset C(S^3/\Gamma)$ is a perimeter minimizing set whose boundary contains $p$. 
By applying Theorem \ref{Theorem | Bernstein}, we can then infer that $C(S^3/\Gamma)$ is isometric to $\bb{R}^4$ and that $F$ is a half space, contradicting $x \in S^E_{N-4} \setminus S^E_{N-5}$. Therefore, $S^E \setminus S^E_{N-5} = \emptyset$ as claimed.

Since $S^E=S^E_{N-5}$, Theorem \ref{Theorem | stratification of a perimeter minimizer} implies that
$
\mathrm{dim}_{\mathcal{H}}(S^E) \leq N-5.
$
\qed
\medskip

Theorem \ref{Theorem| Hausdorff dimension of the singular set P min double bound} is sharp: as shown in the following example, there exists a non-collapsed two-sided Ricci limit space with a perimeter minimizing set $E$ such that $S^E_{N-5}$ is non-empty.

\begin{example} \label{example|Sharpness theorem 1}
    The cone $C(\R\bb{P}^3)$, arising as the blow-down of the Eguchi-Hanson manifold (see \cite[Example 2.15]{CheegerNaberCodim4}),
    is a non-collapsed two-sided Ricci limit space which is singular at the tip. 
    
    By a standard calibration argument it is possible to show that the set 
    \[
    C(\R\bb{P}^3) \times \bb[0,+\infty) \subset C(\R\bb{P}^3) \times \bb{R}
    \]
    is a perimeter minimizing set. Moreover the tip of $C(\R\bb{P}^3)$ belongs to $S^E_{N-5}$.
\end{example}

\printbibliography
\end{document}